\documentclass[12pt]{article}   
\usepackage{amssymb}
\usepackage{latexsym}
\usepackage{color}
\usepackage{amsfonts}
\usepackage{graphicx}
\usepackage{amsthm, amsmath}   

\def\O{{\mathcal O}}

\def\CalC{{\mathcal C}}
\def\P{{\mathbb P}}
\def\G{{\mathbb G}}

\def\Om{\Omega}
\def\qed{\hfill\vbox{\hrule\hbox{\vrule\kern3pt\vbox{\kern6pt}\kern3pt\vrule}\hrule}\bigskip}

\def\qed{\hfill\vbox{\hrule\hbox{\vrule\kern3pt\vbox{\kern6pt}\kern3pt\vrule}\hrule}\bigskip}

\theoremstyle{plain}
\newtheorem{theorem}{Theorem}[section]

\theoremstyle{definition}

\theoremstyle{remark}
\newtheorem{remark}[theorem]{Remark}

\newenvironment{main1}{\bigskip\noindent{\bf Main Result 1}}{\bigskip}
\newenvironment{main2}{\bigskip\noindent{\bf Main Result 2}}{\bigskip}
\newenvironment{pf}{{\it Proof.}}{\qed}
\newenvironment{df}{\bigskip\noindent{\bf Definition}}{\bigskip}

\title{Bertini-type theorems for formal functions in grassmannians}
\author{Jorge Caravantes\\ Universit\'a degli Studi di Genova\\
e-mail: {\tt caravant@dima.unige.it, jorgecaravan@gmail.com}}
\date{}

\begin{document}
\maketitle

\begin{abstract}
The present note generalizes Debarre's Bertini-type results for inverse images of Schubert varieties with the extension of formal functions.
\end{abstract}

\section*{Introduction}

The aim of this short note is to develop Bertini-type results related to the extension of formal functions for inverse images of Schubert varieties of grassmannians. We start with some notation and our fundamental definition.

Given a closed subscheme $Y$ of a variety $X$ over a closed field $k$. Let $X_{/Y}$ be the formal completion of $Y$ along $X$. Let $K(Z)$ be the ring of rational functions of a scheme (formal or not) $Z$ (i.e. the global sections of the sheaf of rational functions). Then it is well known that there exists an injective homomorphism $\alpha:K(X)\to K(X_{/Y})$. We now recall the following definition:

\begin{df}
If $\alpha$ is an isomorphism, we say that $Y$ is G3 in $X$.
\end{df}

The first Bertini-type result regarding extension of formal functions is probably the following:

\begin{theorem}(Faltings, \cite{F})
Let $f:X\to \P^n$ a proper morphism. If dim$(f(X))\ge d$, then for any $d-1$ hyperplanes $H_1,...,H_{d-1}\subset\P^n$, $f^{-1}(H_1\cap...\cap H_{d-1})$ is G3 in $X$.
\end{theorem}

The proof in Faltings' paper involves local rings and local cohomology. Recently, Bonacini, Del Padrone and Nesci found in \cite{BDN} an inspiring global proof whose ideas have already been used for products of projective spaces in \cite{B2}. Our aim is to find an analogous for grassmannians. In order to state our main results, we need a definition regarding special subvarieties of grassmannians.

\begin{df}
Let $\G(k,n)$ be the grassmannian of $k$-planes in $\P^n$. Let $a<b\le n$ be two nonnegative integers. Let $\Pi_0,...,\Pi_k$ be subspaces of $\P^n$ such that dim$\Pi_i=d_i$, with $0\le d_0< ...< d_k\le n$. We introduce the following notation:
\[
\Om(\Pi_0,...,\Pi_k)=\{\Lambda\in\G(k,n)\ |\ \mbox{dim}(X\cap\Pi_i)\ge i,\ \forall i=0,...,k\}
\]
Let also $\Om(d_0,...,d_k)$ be its class in the Chow ring of $\G(k,n)$. They are said to be a Schubert variety and a Schubert cycle (respectively). The dimension of $\Om(d_o,...,d_k)$ is 
\[\sum(d_i-k)=\sum d_i+\frac{k(k+1)}{2}.\]
\end{df}

Now we can state the main results of this note:

\begin{main1} 
{\it Let $X$ be an irreducible complete variety, and let $j\in\{0,...,n\}$. Let $f:X\to \G(k,n)$ be a morphism such that the intersection product $[f(X)]\Om(d_0-1,d_1-1,..,d_j-1,n-k+j+1,...,n)$ in the Chow ring of $\G(k,n)$ is nonzero for some $0<d_0<...<f_j\le n-k+j$. Then for any $\Pi_i\in\G(d_i,n)$, $i=0,...,j$, we have that 
\[f^{-1}(\Om(\Pi_0,...,\Pi_j,\P^{n-k+j+1},...,\P^n))\] 
is G3 in $X$.} 
\end{main1}


The following result is analogous in the case of the grassmannians of lines. It improves the connectedness theorems of Debarre with some complementary conditions (that, in some cases, are strictly milder).

\begin{main2} 
{\it Let $X$ be an irreducible complete variety. Let $f:X\to \G(1,n)$ be a morphism such that the intersection product $[f(X)]\Om(a,a+1)$ in the Chow ring of $\G(1,n)$ is nonzero, given $0\le a<b\le n$. Then for any $\tilde\Lambda\in\G(a,n)$, $\tilde\Pi\in\G(b,n)$, we have $Y:=f^{-1}(\Om(\tilde\Lambda,\tilde\Pi))$ is G3 in $X$.}
\end{main2}


These two Main Results (in the case of grassmannians of lines) can be merged in the following way:

\begin{theorem}
Let $X$ be an irreducible complete variety. Let $f:X\to \G(1,n)$ be a morphism such that one the intersection products 
\begin{itemize}
\item $[f(X)]\Om(a,a+1)$, or
\item $[f(X)]\Om(a-1,b-1)$ (or $[f(X)]\Om(a-1,n)$, if $b=n$)
\end{itemize}
in the Chow ring of $\G(1,n)$ are nonzero, given $0\le a<b\le n$. Then for any $\tilde\Lambda\in\G(a,n)$, $\tilde\Pi\in\G(b,n)$, we have $Y:=f^{-1}(\Om(\tilde\Lambda,\tilde\Pi))$ is G3 in $X$.
\end{theorem}

The paper is distributed as follows. First section containssome results thart are useful for the proofs in the paper. Second section contains the proofs of the main results of this note.

{\bf Acknowledgements:} This paper was written during a two-year stay the Universit\`a degli Studi di Genova with a postdoctoral Fundaci\'on Ram\'on Areces grant. I would like to thank Lucian B\u adescu for introducing me to the problem and Enrique Arrondo for some useful discussions.

\section{Preliminaries}

This section contains several known results that will be needed for the proof of the main result. Some related to G3 property and some related to connectedness.



One of the main tools we can use to prove extension of formal functions is the following theorem:

\begin{theorem}\label{HiroMatsu}(Hironaka-Matsumura, \cite{HM}, see also e.g. \cite[Corollary 9.13]{B1})
Let $f:X'\to X$ be a proper surjective morphism of irreducible varieties over an algebraically closed field $k$. Suppose $Y$ is G3 in $X$. Then $f^{-1}(Y)$ is G3 in $X'$.
\end{theorem}

The following is an useful consequence of Theorem \ref{HiroMatsu}

\begin{theorem}\label{Badescu}(see e.g. \cite[Proposition 9.23]{B1})
Let $f:X'\to X$ be a proper surjective morphism of irreducible varieties over an algebraically closed field $k$, and let $Y\subset X$ and $Y'\subset X'$ be closed subvarieties such that  $f(Y')\subset Y$. Assume that the rings $K(X_{/Y})$ and $K(X'_{/f^{-1}(Y)})$ are both fields. If $Y'$ is G3 in $X'$, then $Y$ is G3 in $X$.
\end{theorem}

In order to use the previous result, we first set a framework where the ring of formal rational functions is a field:

\begin{theorem}\label{connectfield}(see e.g. \cite[Corollary 9.10]{B1})
Let $X$ be an algebraic variety over an algebraically closed field $k$, and let $Y$ be a closed subvariety $X$. Let $u:X'\to X$ be the (birational) normalization of $X$. Then $K(X_{/Y})$ is a field if and only if $u^{-1}(Y)$ is connected.
\end{theorem}

The following result, due to Debarre, is the startup and the inspiration for the two main Bertini-type results, Main Results 1 and 2. It is also used for the proof of the first one (the second having hypotheses milder than this one).

\begin{theorem}(Debarre, \cite[Th\'eor\`eme 8.1 (b)]{D})\label{Debarre}
Let $X$ be a projective variety, $f:X\to \G(k,n)$ a morphism. We consider $0\le d_0<...<d_k\le n$ (we also define $d_{k+1}=n+1$) and fix $\Pi_i\simeq\P^{d_i}\subset\Pi_{i+1}\simeq\P^{d_{i+1}}\subset\P^n$, for all $i=0,...,k-1$. If all the intersection products: 
\begin{itemize}
\item $[f(X)]\Om(d_i-i-1,...,d_{i}-2,d_i-1,d_{i+1},...,d_k)$, for all $i$ such that $d_i>i$ and $d_{i+1}\ne d_i+1$.
\item $[f(X)]\Om(0,...,i,i+1,d_{i+2},...,d_k)$, for all $i$ such that $d_i=i$ and $d_{i+1}\ne i+1$.
\end{itemize} 
in the Chow ring of $\G(1,n)$ are different from zero, then $f^{-1}(\Om(\Pi_0,...,\Pi_k))$ is connected.
\end{theorem}

Again to make Theorems \ref{HiroMatsu} and \ref{Badescu} useful, we need an starting point, which is given by the following result:

\begin{theorem}\label{G3P}(Hironaca-Matsumura \cite{HM} for $\P^n$, Babakarian-Hironaka \cite{BH} for general case)
Any positive-dimensional variety in a grassmannian (not necessarily of lines) is G3.
\end{theorem}

\section{Proofs of the results}

In this section we prove the two main results of this note, just after sketching the general method.

\begin{remark}\label{metodo}
We introduce in this remark the general method we use to prove Bertini-type results regarding extension of formal functions, as extracted from \cite{BDN} and \cite{B2}. 

We start with a closed subscheme $Y$ of a scheme $X$. In order to prove that $Y$ is G3 in $X$, we consider an incidence variety $Z$ and the incidence diagram:
\[\begin{array}{ccccl}
& &Z&\subset &X\times P\\
&\swarrow_p & & \searrow^q\\
X& & & & P
\end{array}\]
with $P$ a sufficiently good space (in our cases a projective space and a grassmannian).

We sketch the steps we follow for the proofs of the main results:
\begin{enumerate}
\item prove that $Z$ is irreducible.
\item prove that both $p$ and $q$ are proper and surjective.
\item find a subscheme $L\subset P$ satisfying:
\begin{itemize}
\item $L$ is G3 in $P$.
\item $Y':=q^{-1}(L)\subset p^{-1}(Y)$.
\end{itemize}
\item when the previous steps are completed, we automatically get that $Y'$ is G3 in $Z$ by Theorem \ref{HiroMatsu} to the projection $q:Z\to P$. 
\item prove that both $K(X_{/Y})$ and $K(Z_{/Y'})$ are fields.
\item to finish, by Theorem \ref{Badescu} applied to the projection $p:Z\to X$, we get that $Y$ is G3 in $X$. 
\end{enumerate}
\end{remark}

We finish the section with the proofs of the main results.


\begin{pf} \emph{(of Main Result 1)}
We will write the proof for the case $j=k$, the other being similar.

Using Theorem \ref{HiroMatsu}, we can suppose $X\subset \G(k,n)$ and $f$ to be the inclusion. Let us define, as in Remark \ref{metodo}, the incidence variety:
\[\begin{array}{ccccl}
&&Z& := & \left\{ (\Lambda,[s_{d_0},...,s_n])\ \left|\ 
{\scriptsize\begin{array}{l}
s_{d_k|\Lambda}=...=s_{n|\Lambda}=0\\ rk(s_{d_i|l},...,s_{d_{i+1}-1|l})\le 1,\mbox{ when } d_i\ne d_{i+1}-1 
\end{array}}\right. \right\} \\
& \swarrow_p && \searrow^q\\
X& & & & P=\P(H^0(\P^n,\O_{\P^n}(1)))
\end{array}\]
Obviously, when $d_i=d_{i+1}-1$, the rank 1 condition means nothing (moreover, in that case, $\Pi_i$ does not impose condition for $\Lambda$). In the case $j<n$ we would just omit the vanishing condition.

First step is to prove that $Z$ is irreducible. We first need to describe $Z$ better. The fiber of $p$ on a point $\Lambda\in X$ consists of all $(n-d_0+1)$-tuples of linear sections such that, when restricted to the line $l$, the last $n-d_k+1$ sections vanish and any pack $s_{d_i},...,s_{d_{i+1}-1}$ have rank at most one, all after restriciton to $\Lambda$. This makes $p$ a bundle over $X$. To study the fiber, we put $x_0,...,x_n$ as coordinates in $\P^n$ such that the space $\Lambda_0:=V(x_{k+1},...,x_n)$ lies in $X$. The fiber of $p$ in $\Lambda_0$ consists on packs of $n-d_0+1$ sections that are linear combinations of the variables $x_0,...,x_n$. The coefficients of $x_k+1,...,x_n$ are completely free, but the coefficients of $x_0,...,x_k$:
\begin{itemize}
\item vanish for $s_{d_k},...,s_n$.
\item form a rank 1 matrix for any pack $s_{d_i},...,s_{d_{i+1}-1}$.
\end{itemize} 
This means that the ideal of $p^{-1}(\Lambda_0)$ in the graded ring of $P$ coincides with the ideal of
\[\CalC(\P^k\times\P^{d_2-d_1-1})\times...\times\CalC(\P^k\times\P^{d_k-d_{k-1}-1})\times\P^{(n-k)(n-d_k+1)-1}\]
(where $\CalC(Y)$ is the cone of a variety $Y$) in 
\[\P^{(d_2-d_1)(k+1)-1}\times...\times\P^{(d_k-d_{k-1})(k+1)-1}\times\P^{(n-d_k+1)(k+1)-1}\]
(which has the same ring of $P$, but is differently graded). Such ideal is irreducible. Then $Z$ is a bundle with irreducible fibers and then irreducible.



To complete second step we must show that both $p$ and $q$ are surjective. The first one defines a bundle on $X$, so it is obviously surjective. The fiber of $q$ at a general point $[t_{d_0},...,t_n]$, $t_i\in H^0(\O_{\P^n}(1))$, is clearly the set of points in $X$ representing a $k$-plane that is contained in $V(t_{d_k+1},...,t_n)\simeq\P^{d_k-1}$ and intersects all $V(t_{d_i+1},...,t_n)\simeq\P^{d_i-1}$ in a $\P^i$. This fiber is not empty due to the condition: $[f(X)]\Om(d_0-1,...,d_k-1)\ne0$.

Step 3 consists of defining the subscheme $L$ in $P$ with its two properties. Let $t_{a+1},...,t_n\in H^0(\P^n,\O_{\P^n}(1))$ satisfy $\Pi_i=V(t_{d_i+1},...,t_n)$ for all $i=0,...,k$. Then we consider the line $L=<(0:t_{d_0+1}:...:t_n),(t_{d_0+1}:...:t_n:0)>$. For any $(s_{d_0}:...:s_n)\in N$, we have that $V(s_{d_i},...,s_n)\subset\tilde\Pi_i$, so $q^{-1}(N)\subset p^{-1}(Y)$. On the other side:
\[q^{-1}(t_{a+1}:...:t_n:0)=\{(t_{a+1}:...:t_n:0)\}\times Y\]
since $V(t_{d_{i}+1},...,t_n,0)=\Pi_i$. Therefore $p(q^{-1}(N))= Y$. Since, by Theorem \ref{G3P}, $N$ is G3 in a projective space, this proves the connectedness of $Y$.

Now, we directly get fourth step, for $Y':=q^{-1}(N)$ is G3 in $Z$ by Theorem \ref{HiroMatsu}.

We continue proving that both $K(Z_{/p^{-1}(Y)})$ and $K(X_{/Y})$ are fields. Let $h:\tilde Z\to Z$ and $g:\tilde X\to X$ birational normalizaions of $Z$ and $X$ respectively. Then, by Theorem \ref{Debarre}, both $h^{-1}p^{-1}(Y)$ and $g^{-1}(Y)$ are connected, since they are preimages of $\Om(\Pi_0,...,\Pi_k)$. Now, by Theorem \ref{connectfield} $K(Z_{/p^{-1}(Y)})$ and $K(X_{/Y})$ are fields

Last step is again automatic, since by Theorem \ref{Badescu}, $Y'$ being G3 in $Z$ implies that $Y$ is G3 in $X$.
\end{pf}

\medskip

We prove now Main Result 2.



\begin{pf} \emph{(of Main Result 2)}
First of all, we can take a resolution of the singularities of $X$ and compose with $f$. Applying Theorem \ref{HiroMatsu}, we can suppose $X$ to be smooth. We again follow the program in Remark \ref{metodo} and define the incidence variety
$$\begin{array}{ccccl}
&&Z& := & \{(x,\Gamma)\ |\ f(x)=L\subset\Gamma\} \\
& \swarrow_p && \searrow^q\\
X& & & & P=\G(a+1,n)
\end{array}$$

To begin the program we show that $Z$ is irreducible. The fiber of $p$ on $L$ is the set of all $\Gamma$ that contain the line $L=f(x)$ which is isomorphic to a $\G(a,n-2)$. Therefore $Z$ is a bundle over $X$ with irreducible and smooth fiber. We so get that $Z$ is both irreducible and smooth.

We have seen moreover that $p$ is surjective. On the other side, the inverse image of a given $\Gamma$ by $q$ is the preimage of the set of all lines $L\in f(X)$ contained in $\Gamma$. Since $[f(X)]$ intersects $\Om(a,a+1)$, it is nonempty and hence $q$ is also surjective. This finishes the second step. 

For the third step of the program, we consider the subspace $N\subset\G(a+1,n)$ consisting of all $\Gamma$ such that $\tilde\Lambda \subset \Gamma \subset \tilde\Pi$. Clearly, $N\simeq\P^{b-a-1}$, so $N$ is G3 in $\G(a+1,n)$ by Theorem \ref{G3P}.

To finish this step, we need to prove that $q^{-1}(N)\subset p^{-1}(Y)$. For any subspace $\Gamma\in L$ we have that its fiber by $q$ consists on all points $x\in X$ such that $f(x)=L\subset\Gamma$. We know that $\Gamma\subset\tilde\Pi$ (by definition of $N$). Moreover, since $\tilde\Lambda$ is a hyperplane of $\Gamma$ (also by definition of $N$), we have that every line of $\Gamma$ intersects $\tilde\Lambda$, so all $x\in p(q^{-1}(\{\Gamma\}))$ satisfy $L=f(x)\subset\tilde\Pi$ and $L\cap\tilde\Lambda\ne\emptyset$, so $p(q^{-1}(\{\Gamma\}))\subset Y$ as we wanted.

Fourth step is automatic, $q^{-1}(N)$ is G3 in $Z$ by Theorem \ref{HiroMatsu} and the fact that $N$ is G3 in $\G(a+1,n)$.


Next step is the proof of both $K(Z_{/p^{-1}(Y)})$ and $K(X_{/Y})$ being fields. Since both $X$ and $Z$ are smooth (and then normal), Theorem \ref{connectfield} reduces this question to showing that both $Y$ and $p^{-1}(Y)$ are connected. Moreover, $p$ is a bundle, so proving that $Y$ is connected is enough: For all $x\in Y$, $l=f(x)$ satisfies $L\cap\tilde\Lambda\ne\emptyset$, so $\Gamma:=<L,\tilde\Lambda>$ has dimension $b-1$ and therefore $L\in p(q^{-1}(\{\Gamma\}))$. This proves that $Y\subset p(q^{-1}(N))$. Since the other content was proved in the third step, we have $Y=p(q^{-1}(N)=Y')$. $Y'$ is G3 in $Z$ as shown before, so it is conected. Therefore $Y$ is the image of a connected subset by a continuous function, so it is also connected.

Last step is automatic once more, since by Theorem \ref{Badescu},  $Y'$ is G3 in $Z$ provides the fact that $Y$ is G3 in $X$.
\end{pf}


\begin{remark}
The generalization of Main Result 2 to general grassmannians requires better knowledge on Flag manifolds than that of the author. So an unique statement for general grassmannians is out of reach now.
\end{remark}

\end{document}